\documentclass[12pt, oneside]{article}
\usepackage{amsthm}
\usepackage{graphicx} 
\usepackage{tikz}
\usetikzlibrary{calc,intersections,through} 
\usetikzlibrary{decorations.fractals}
\usepackage[english]{babel}
\usepackage[utf8]{inputenc}
\usepackage[many]{tcolorbox}
\usepackage[a4paper, total={6.5in,10in}]{geometry}
\usepackage{amsfonts,amsmath}
\usepackage{amsmath}
\usepackage{tkz-euclide}
\usepackage{float}
\usepackage{relsize}
\newtheorem{definition}{Definition}
\newtheorem{remark}{Remark}
\newtheorem{proposition}{Proposition}
\AtBeginDocument{}
\usepackage{setspace}
\newtheorem{theorem}{Theorem}[section]
\newtheorem{lemma}[theorem]{Lemma}
\usepackage{mathtools}

\title{Riemann Integration in $\mathbb{R}^n$}
\author{Akerele Olofin Segun\\ Department of Mathematics, University of Ibadan, Ibadan\\ sakerele647@stu.ui.edu.ng, akereleolofin@gmail.com}

\begin{document}
	\maketitle
	\begin{abstract}
	\par During this quest, we define the Riemann integrals using the Darboux upper and lower integrals. The ideas here are very similar to integration in one dimension. The complication is mostly notational. \\
	The differences between one and several dimensions will grow more pronounced as we proceed.
	\par The so-called Riemann sums have their origin in the efforts of Greek mathematicians to find the center of gravity or the volume of a solid body. These researches led to the method of exhaustion, discovered by Archimedes and described using modern ideas by MacLaurin in his \textit{Treatise of Fluxions} in 1742. At this times the sums were only a practical method for computing an area under a curve, and the existence of this area was considered geometrically obvious. The method of exhaustion consists in almost covering the space enclosed by the curve with $n$ geometric objects with well-known areas such as rectangles or triangles, and finding the limit (though this topic was very blurry at these early times) when $n$ increases. One of its most remarkable application is squaring the area $\mathcal{A}$ enclosed by a parabola and a line.
	 Nowadays, Riemann sums remain a useful tool to study some sequences involving sums.
	 \footnote{Bernhard Riemann : \emph{The greatest strategy is doomed if it's implemented badly.}}
\end{abstract}
	\newpage	
	\section*{Preliminaries}
	\section{The Riemann Integral in $\mathbb{R}$}
	\subsection{Definition.}
	Let $[a,b]$ be a given closed and bounded interval in $\mathbb{R}$. A partition $\mathcal{P}$ of $[a,b]$ is a finite set of points $\mathcal{P}=\{x_0,x_1,x_2,\dots,x_n\}$ such that 
	$$ a=x_0<x_1<x_2<\cdots < x_n=b $$
	\subsection{Remark.}
	\begin{itemize}
		\item [1.] There is no requirement that the partition points $x_i$ be equally spaced. Partitions in which the partition points are equally spaced are called standard partitions.
		\item[2.] For each $i=1,2,3,\dots,n$. Set 
		$$\Delta x_i = x_i - x_{i-1}$$
		which is equal to the length of the subinterval $[x_{i-1},x_i]$.
		\item [3.] The number $\|\mathcal{P}\|=\max_{1\leq i \leq n}\Delta x_i$ is called the norm of $\mathcal{P}$ (or the mesh of $\mathcal{P}$).
		\item [4.] Suppose that $f$ is a bounded function on $[a,b]$, we write 
		\begin{align*}
			M_i(f)=\sup\{f(x):x\in [x_{i-1},x_i]\}\\
			m_i(f)=\inf\{f(x):x\in [x_{i-1},x_i]\}
		\end{align*}
	\end{itemize}
	\subsection{Definition.}
	Let $f$ be a bounded function on $[a,b]$ and $\mathcal{P}=\{x_0,x_1,x_2,\dots,x_n\}$ a partition of $[a,b]$. The upper Riemann sum for $f$ and the partition $\mathcal{P}$ is defined by 
	$$U(f,\mathcal{P})=\sum_{i=1}^{n}M_i(f)\Delta x_i $$
	Similarly, the lower Riemann sum for $f$ and the partition $\mathcal{P}$ is defined by 
	$$L(f,\mathcal{P})=\sum_{i=1}^{n}m_i(f)\Delta x_i $$
	\subsection*{Remark.}	
	\begin{itemize}
		\item [1.] Since $m_i(f)\leq M_i(f)$ for all $i=1,2,3,\dots, n$, we always have 
		$$L(f,\mathcal{P})\leq U(f,\mathcal{P}).$$
		for any partition $\mathcal{P}$ of $[a,b]$
		\item [2.] For a non-negative continuous function, the upper Riemann sum $U(f,\mathcal{P})$ represents the circumscribed rectangular approximation to the area under the graph $f$. Similarly, the lower Riemann sum represents the inscribed rectangular approximation of the area under the graph of $f$.
	\end{itemize}

	\begin{lemma}\label{lm2}
		Let $f$ be a bounded function on $[a,b]$. If $m\leq f(x)\leq M$ for all $x\in[a,b]$, then 
		$$m(b-a)\leq L(f,\mathcal{P}) \leq U(f,\mathcal{P})\leq M(b-a)$$
		for any partition $\mathcal{P}$ of $[a,b]$
	\end{lemma}
	\begin{proof}
		Let $\mathcal{P}=\{x_0,x_1,x_2,\dots,x_n\}$ be any partition of $[a,b]$. Since $M_i(f)\leq M$ for all $i=1,2,3,\dots,n$, we have
		$$U(f,\mathcal{P})=\sum_{i=1}^{n}M_i(f)\Delta x_i \leq \sum_{i=1}^n M(x_i-x_{i-1})=M(b-a)$$
		Similarly, since $m_i(f)\geq m$ for all $i=1,2,3,\dots,n$, we have 
		$$L(f,\mathcal{P})=\sum_{i=1}^{n}m_i(f)\Delta x_i \geq \sum_{i=1}^n m(x_i-x_{i-1})=m(b-a)$$
		Hence, 
		$$m(b-a)\leq L(f,\mathcal{P}) \leq U(f,\mathcal{P})\leq M(b-a)$$
		for any partition $\mathcal{P}$ of $[a,b]$
	\end{proof}
	\subsection*{Definition}	
	A partition $\mathcal{P}'$ of $[a,b]$ is called a refinement of $\mathcal{P}$ if $\mathcal{P}\subset \mathcal{P}'$.
	\begin{lemma}\label{lm1}
		Let $f$ be bounded on $[a,b]$ and let $\mathcal{P}'$ be a refinement of a partition $\mathcal{P}$ of $[a,b]$. Then
		$$L(f,\mathcal{P})\leq L(f,\mathcal{P}')\leq U(f,\mathcal{P}')\leq U(f,\mathcal{P})$$
	\end{lemma}
	\begin{proof}
		Let $\mathcal{P}=\{x_0,x_1,x_2,\dots,x_n\}$. Suppose that $\mathcal{P}'=\mathcal{P}\cup \{x'\}$. Furthermore, we can, without the loss of generality, assume that $x_0<x'<x_1$, so that $\mathcal{P}'=\{x_0,x',x_1,x_2,\dots,x_n\}$.\\
		Now, let 
		\begin{align*}
			M'_1(f)=\sup\{f(x):x\in [x_0,x']\}\\
			M''_1(f)=\sup\{f(x):x\in [x',x_1]\}	
		\end{align*}
		Since $f(x)\leq M_1(f)$ for all $x\in [x_0,x_1]$, we have that $f(x)\leq M_1(f)$ for all $x\in [x_0,x']$ and also for all $x\in [x',x_1]$. Therefore,
		$$M'_1(f)\leq M_1(f) \quad and \quad M''_1(f)\leq M_1(f)$$
		Hence,
		$$M'_1(f)(x'-x_0)+M''_1(f)(x_1-x')\leq M_1(f)(x'-x_0)+M_1(f)(x_1-x') = M_1(f)(x_1-x_0)$$
		Now,
		\begin{align*}
			U(f,\mathcal{P}')=M'_1(f)(x'-x_0)+M''_1(f)(x_1-x')+\sum_{i=2}^nM_i(f)\Delta x_i \leq M_1(f)(x_1-x_0)+\sum_{i=2}^{n}M_i(f)\Delta x_i\\
			=\sum_{i=1}^nM_i(f)\Delta x_i = U(f,\mathcal{P}).
		\end{align*}
		Showing that $U(f,\mathcal{P}')\leq U(f,\mathcal{P})$.\\
		If $\mathcal{P}'$ contains more than one additional point, we repeat the above argument the appropriate number of times. That $L(f,\mathcal{P})\leq L(f,\mathcal{P}')$ is proved in a similar way.
	\end{proof}
	\subsection*{Definition}
	Let $f$ be a bounded real-valued function on a closed and bounded interval $[a,b]$. The upper and lower Riemann integrals of $f$, denoted by 
	$$\overline{\int_a^b}f(x)\,dx \quad and \quad \underline{\int_a^b}f(x)\,dx$$
	respectively are defined by 
	\begin{align*}
		\overline{\int_a^b}f(x)\,dx=\inf\{U(f,\mathcal{P}):\mathcal{P}\in \mathbb{P}[a,b]\}\\
		\\
		\underline{\int_a^b}f(x)\,dx = \sup\{L(f,\mathcal{P}):\mathcal{P}\in \mathbb{P}[a,b]\}
	\end{align*}
	\begin{theorem}
		Let $f$ be a bounded real-valued function on $[a,b]$. Then 
		$$ \underline{\int_a^b}f(x)\,dx\leq \overline{\int_a^b}f(x)\,dx $$
	\end{theorem}
	\begin{proof}
		Let $P_1$ and $P_2$ be any two partitions of $[a,b]$. Then from Lemma \ref{lm1}
		$$L(f,P_1)\leq L(f,P_1\cup P_2)\leq U(f,P_1\cup P_2) \leq U(f,P_2)$$
		Thus,
		$$L(f,P_1)\leq U(f,P_2)$$
		for any two partitions $P_1$ and $P_2$. Hence, 
		$$\underline{\int_a^b}f(x)\,dx = \sup_{P_1}L(f,P_1)\leq U(f,P_2).$$
		That is 
		$$\underline{\int_a^b}f(x)\,dx = U(f,P_2).$$
		Taking the infimum over $P_2$, we obtain 
		$$\underline{\int_a^b}f(x)\,dx \leq \inf_{P_2}U(f,P_2)=	\overline{\int_a^b}f(x)\,dx.$$
		Hence, the result.
	\end{proof}
	\subsection*{Definition}
	Let $f$ be a bounded real-valued function on a closed and bounded interval $[a,b]$. $f$ is called \textbf{\textit{Riemann-Integrable}} on $[a,b]$ if 
	$$\underline{\int_a^b}f(x)\,dx=\overline{\int_a^b}f(x)\,dx.$$
	The common value is denoted by, 
	$$\int_a^b f(x)\, dx$$	
	and is called the \textit{Riemann integral} of $f$ on $[a,b]$.\\
	\par We denote by $\mathcal{R}[a,b]$ the set of all Riemann-integrable functions on $[a,b]$. 
	\subsection{Remark.}
	If $f:[a,b]\longrightarrow \mathbb{R}$ satisfies $m\leq f(x) \leq M$ for all $x\in [a,b]$, then by Lemma \ref{lm2}, 
	$$m(b-a)\leq \underline{\int_a^b}f(x)\,dx \leq \overline{\int_a^b}f(x)\,dx \leq M(b-a).$$
	If in addition $f\in \mathcal{R}[a,b]$, then 
	$$m(b-a)\leq \int_a^b f(x)\, dx\leq M(b-a).$$
	In particular, if $f(x)\geq 0$ for all $x\in [a,b]$, then $\displaystyle{\int_a^b f(x)\,dx \geq 0}$. If $f\in \mathcal{R}[a,b]$ is non-negative, then the quantity  $\displaystyle{\int_a^b f(x)\,dx \geq 0}$ represents the area of the region bounded above by the graph $y=f(x)$, below by the $x-axis$ and by the lines $x=a$ and $x=b$.	
	\section{Riemann Integral Over Rectangles}
	
	\subsection{Rectangles and Partitions}
	\subsubsection*{Definition.}
	Let $(\xi_1,\xi_2,\dots,\xi_n)$ and $(\eta_1,\eta_2,\dots,\eta_n)$ be such that $\xi_k \leq \eta_k$ for all $k$. A set of the form $[\xi_1,\eta_1]\times [\xi_2,\eta_2]\times \cdots \times [\xi_n,\eta_n]$ is called a closed rectangle. In this setting it is sometimes useful to allow $\xi_k=\eta_k$, in which case we think of $[\xi_k,\eta_k]=\{\xi_k\}$ as usual. If $\xi_k < \eta_k$ for all $k$, then a set of the form $(\xi_1,\eta_1)\times (\xi_2,\eta_2)\times \cdots \times (\xi_n,\eta_n)$ is called an open rectangle. \\
	For any closed or open rectangle $\textbf{R}:= [\xi_1,\eta_1]\times [\xi_2,\eta_2]\times \cdots \times [\xi_n,\eta_n]\subset \mathbb{R}^n$ or\\ $\textbf{R}:= (\xi_1,\eta_1)\times (\xi_2,\eta_2)\times \cdots \times (\xi_n,\eta_n)\subset \mathbb{R}^n $, we define the \textit{n-dimensional volume} by 
	\begin{equation}
		V(\textbf{R}):=(\eta_1-\xi_1)(\eta_2-\xi_2)\cdots (\eta_n-\xi_n)=\prod_{i=1}^n(\eta_i-\xi_i)
	\end{equation}
	\par	A partition $P$ of the closed rectangle $\textbf{R}=[\xi_1,\eta_1]\times [\xi_2,\eta_2]\times \cdots \times [\xi_n,\eta_n]$ is a finite set of partitions $P_1,P_2,\dots,P_n$ of the intervals $[\xi_1,\eta_1], [\xi_2,\eta_2], \dots, [\xi_n,\eta_n]$. We will write $P=(P_1,P_2,\dots,P_n)$. That is, for every $k$ there is an integer $\ell_k$ and the finite set of numbers $P_k=\{x_{k,0},x_{k,1},x_{k,2},\dots,x_{k,\ell_k}\}$ such that 
	$$\xi_k = x_{k,0}<x_{k,2}<\cdots < x_{k,\ell_{k-1}}<x_{k,\ell_k}=\eta_k$$
	Picking a set of $n$ integers $j_1,j_2,\dots,j_n$ where $j_k\in \{1,2,3,\dots,\ell_k\}$ wee get the sub-rectangle
	$$[x_{1,j_1-1},x_{1,j_1}]\times [x_{2,j_2-1},x_{2,j_2}] \times \cdots \times [x_{n,j_n-1},x_{n,j_n}]$$
	For simplicity, we order the sub-rectangles somehow and we say $\{R_1,R_2,R_3,\dots,R_N\}$	are the sub-rectangles corresponding to the partition $P$ of $R$. It is not difficult to see that these sub-rectangles cover our original $R$, and their volume sums to that of $R$. That is 
	\begin{align}
		R=\bigcup_{j=1}^NR_j, \quad and \quad V(R)=\sum_{j=1}^{N}V(R_j)
	\end{align}
	When 
	$$R_k=	[x_{1,j_1-1},x_{1,j_1}]\times [x_{2,j_2-1},x_{2,j_2}] \times \cdots \times [x_{n,j_n-1},x_{n,j_n}]$$
	then
	$$V(R_k)=\Delta x_{1,j_1}\Delta x_{2,j_2}\cdots \Delta x_{n,j_n}=\prod_{i=1}^n\Delta x_{i,j_i}=\prod_{i=1}^n(x_{i,j_i}-x_{i,j_{i}-1}).  $$
	\par Let $R\subset \mathbb{R}^n$ be a closed rectangle and let $f:R\longrightarrow \mathbb{R}$ be a bounded function, Let $P$ be partition of $[a,b]$ and suppose that there are $N$ sub-rectangles. Let $R_i$ be a sub-rectangle of $P$. Define
	
	\begin{align*}
		m_i := \inf\{f(x):x\in R_i\},\\
		M_i := \sup\{f(x):x\in R_i\},\\
		L(P,f):= \sum_{i=1}^N m_iV(R_i),\\
		U(P,f):= \sum_{i=1}^N M_iV(R_i).
	\end{align*}
	We call $L(P,f)$ the lower Darboux sum and $U(P,f)$ the upper Darboux sum.
	\par The indexing in the definition may be complicated, fortunately we generally do not need to go back directly to the definition often. We start proving facts about the Darboux sums analogous to the one-variable results.\\
	
	\begin{proposition}\label{lm3}
		Suppose $R\subset\mathbb{R}^n$ is a closed rectangle and $f:R\longrightarrow \mathbb{R}$ is a bounded function. Let $m,M \in \mathbb{R}$ be such that for all $\xi \in R$ we have $m\leq f(\xi) \leq M$. For any partition P of R we have
		$$mV(R)\leq L(P,f) \leq U(P,f)\leq MV(R).$$
	\end{proposition}
	\begin{proof}
		Let $P$ be a partition. Then note that $m\leq m_i$ for all i. Also $m\leq M_i$ for all i. finally $\sum_{i=1}^N V(R_i)=V(R).$\\ Therefore,\\
		\begin{align*}mV(R)=m\left(\sum_{i=1}^N V(R_i)\right)=\sum_{i=1}^N mV(R_i) \leq \sum_{i=1}^N m_iV(R_i) \leq \sum_{i=1}^N M_iV(Ri) \leq \sum_{i=1}^N MV(R_i) \\ =M\left( \sum_{i=1}^N V(R_i)\right)=MV(R)
		\end{align*}
		
	\end{proof}
	\section{Upper and Lower Integrals}
	By proposition \ref{lm3} the set of upper and lower Darboux sums are bounded sets and we can take their infima and suprema. As before, we now make the following definition.
	\begin{definition}
		If $f:R\longrightarrow \mathbb{R}$ is a bounded function an a closed rectangle $R \subset \mathbb{R}^n$. Define
		$$\underline{\int_{R}}f:= \sup\{L(P,f): \text{P a partition of R}\}, \qquad \overline{\int_{R}}f:= \inf\{U(P,f): \text{P a partition of R}\}, $$
		We call $\displaystyle{\underline{\int}}$ the lower Darboux integral and $\displaystyle{\overline{\int}}$ the upper Darboux integral.
		\par As in one dimension we have refinements of partitions.
		
	\end{definition}
	
	\begin{definition}
		Let $R \subset \mathbb{R}^n$ be a closed rectangle and let $P=(P_1,P_2,\dots,P_n)$ and $\tilde{P}=(\tilde{P}_1,\tilde{P}_2,\dots,\tilde{P}_n) $ be partitions of R. We say $\tilde{P}$ a refinement of P if as sets $P_k \subset \tilde{P}_k$ for all $k=1,2,\dots,n$.
	\end{definition}
	\par It is not difficult to see that if $\tilde{P}$ is a refinement of P, then sub-rectangles of P are unions of sub-rectangles of $\tilde{P}$. Simply put, in a refinement we took the sub-rectangles of P and we cut them into smaller sub-rectangles.
	\begin{proposition}\label{p2}
		Suppose $R\subset \mathbb{R}^n$ is a closed rectangle, $P$ is a partition of $R$ and $\tilde{P}$ is a refinement of $P$. If $f:R\longrightarrow \mathbb{R}$ be a bounded function, then 
		$$L(P.f)\leq L(\tilde{P},f) \quad and \quad U(\tilde{P},f)\leq U(P,f)$$
	\end{proposition}
	\begin{proof}
		Let $R_1,R_2,R_3,\dots, R_N$ be the sub-rectangles of $P$ and $\tilde{R_1},\tilde{R_2},\dots, \tilde{R}_M$ be the sub-rectangles of $\tilde{R}$. Let $I_k$ be the set of indices $j$ such that $\tilde{R}_j \subset R_k$. We notice that 
		$$R_k = \bigcup_{j\in I_k}\tilde{R}_j, \qquad V(R_k)=\sum_{j\in I_k}V(\tilde{R}_j).$$
		Let $m_j:=\inf\{f(\xi): \xi \in R_j\}$, and $\tilde{m}_j:=\inf\{f(\xi): \xi \in \tilde{R}_j\}$ as usual. Notice also that if $j\in I_k$, then $m_k\leq \tilde{m}_j$. Then 
		$$L(P,f)=\sum_{k=1}^Nm_kV(R_k)=\sum_{k=1}^N\sum_{j\in I_k}m_kV(\tilde{R}_j)\leq \sum_{k=1}^N\sum_{j\in I_k}\tilde{m}_jV(\tilde{R}_j)=L(\tilde{P},f).$$
	\end{proof}	
	The key point of this next proposition is that the lower Darboux integral is less than or equal to the upper Darboux integral.
	\begin{proposition}\label{p3}
		Let $R\subset \mathbb{R}^n$ be a closed rectangle and $f:R\longrightarrow \mathbb{R}$ a bounded function. Let $m,M \in \mathbb{R}$ be such that for all $\xi \in R$ we have $m\leq f(\xi)\leq M$. Then
		\begin{equation}
			mV(R)\leq \underline{\int_R}f \leq \overline{\int_R}f \leq MV(R)
		\end{equation}
	\end{proposition}
	\begin{proof}
		Let $P$ be a partition, via proposition \ref{lm3} 
		$$mV(R)\leq L(P,f)\leq U(P,f)\leq MV(R).$$
		By taking suprema of $L(P,f)$ and infima of $U(P,f)$ over all $P$ we obtain the first and the last inequality. Now, let $P=(P_1,P_2,\dots,P_n)$ and $Q=(Q_1,Q_2,\dots,Q_n)$ be partitions of $R$. Define $\tilde{P}=(\tilde{P}_1.\tilde{P}_2,\dots,\tilde{P}_n)$ by letting $\tilde{P}_k = P_k \cup Q_k$. Then $\tilde{P}$ is a partition of $R$ as can easily be checked, and $\tilde{P}$ is a refinement of $P$ and a refinement of $Q$. By proposition \ref{p2}, $L(P,f)\leq L(\tilde{P},f)$ and $U(\tilde{P},f)\leq U(Q,f)$. Therefore,
		$$L(P,f)\leq L(\tilde{P},f)\leq U(\tilde{P},f)\leq U(Q,f).$$
		It follows directly that,
		\vspace{-4mm}
		$$\sup\{L(P,f): \text{P a partition of R}\} \leq  \inf\{U(P,f): \text{P a partition of R}\}.$$
		In other words,
		$$	mV(R)\leq \underline{\int_R}f \leq \overline{\int_R}f \leq MV(R)$$
	\end{proof}		
	\section{The Riemann Integral in $\mathbb{R}^n$}
	\par 	We now have all we need to define the Riemann integral in n-dimensions over rectangles. Again, the Riemann integral is only defined on a certain class of functions, called the Riemann integrable functions.
	\begin{definition}
		Let $R\subset \mathbb{R}^n$ be a closed rectangle. Let $f:R\longrightarrow \mathbb{R}$ be a bounded function such that
		$$ \underline{\int_R}f(x)\, dx = \overline{\int_R}f(x)\,dx.$$
		Then $f$ is said to be Riemann integrable. The set of Riemann integrable functions on $R$ is denoted by $\mathcal{R}(R)$. When $f\in \mathcal{R}(R)$ we define the Riemann integral 
		$$ \int_R f := \underline{\int_R}f=\overline{\int_R}f.$$
		When the variable $\xi \in \mathbb{R}^n$ needs to be emphasized we write 
		$$\int_Rf(\xi)\,d\xi, \quad \int_Rf(x_1,x_2,x_3,\dots,x_n)\,dx_1\cdots dx_n, \quad or \quad \int_R f(\xi)\,dV.$$
		If $R\subset \mathbb{R}^2$, then often instead of volume we say area, and hence write 
		$$\int_Rf(\xi)\,dA.$$
	\end{definition}
	\begin{remark}
		Let $f:R\longrightarrow \mathbb{R}$ be a Riemann integrable function on a closed rectangle $R\subset \mathbb{R}^n$. Let $m,M\in \mathbb{R}$ be such that $m\leq f(\xi)\leq M$ for all $\xi \in R$. Then
		$$mV(R)\leq \int_R f \leq MV(R).$$
	\end{remark}
	\par For a example a constant function is Riemann integrable. Suppose $f(\xi)=k$ for all $\xi \in R$. Then 
	$$kV(R)\leq \underline{\int_R}f\leq \overline{\int_R}f \leq kV(R).$$	
	So $f$ is integrable, and futhermore 
	$$\int_R f =cV(R).$$
	\begin{remark}(Linearity).
		Let $R\subset \mathbb{R}^n$ be a closed rectangle and let $f$ and $g$ be in $\mathcal{R}(R)$ and $\alpha \in \mathbb{R}$.
		\begin{itemize}
			\item [1.] $\alpha f$ is in $\mathcal{R}(R)$ and 
			$$\int_R \alpha f = \alpha \int_R f$$
			\item [2.] $f+g$ is in $\mathcal{R}(R)$ and 
			$$\int_R(f+g)=\int_R f + \int_R g.$$
		\end{itemize}
	\end{remark}
	\begin{remark}(Monotonicity).
		Let $R\subset \mathbb{R}^n$ be a closed rectangle and let $f$ and $g$ be in $\mathcal{R}(R)$ and let $f(\xi)\leq g(\xi)$ for all $\xi \in R$. Then
		$$\int_R f\leq \int_Rf|_R.$$
	\end{remark}	
	The proofs of linearity and monotonicity are almost completely identical as the proofs from one variable. We therefore omit the proofs.	
	\begin{proposition}
		For a closed rectangle $S\subset \mathbb{R}^n$, if $f:S\longrightarrow \mathbb{R}$ is integrable and $R\subset S$ is a closed rectangle, then $f$ is integrable over $R$.
	\end{proposition}	
	\begin{proof}
		Given $\varepsilon >0$, we find a partition $P$ such that $U(P,f)-L(P,f)<\varepsilon$. By making a refinement of $P$ we can assume that the endpoints of $R$ are in $P$, or in other words, $R$ is a union of sub-rectangles of $P$. Then the sub-rectangles of $P$ divide into two collections, ones that are subsets of $R$ and ones whose intersection with the interior of $R$ is empty. Suppose that $R_1,R_2,\dots,R_K$ be the sub-rectangles that are subsets of $R$ and $R_{K+1},\dots,R_N$ be the rest. Let $\tilde{P}$ be the partition of $R$ composed of those sub-rectangles of $P$ contained in $R$. Then using the same notation as before.
		\begin{align*}
			\varepsilon > U(P,f)-L(P,f)=\sum_{k=1}^K(M_k-m_k)V(R_k)+\sum_{k=K+1}^{N}(M_k-m_k)V(R_k)\\
			\geq \sum_{k=1}^K(M_k-m_k)V(R_k)=U(\tilde{P},f|_R)-L(\tilde{P},f|_R)
		\end{align*}
		Therefore $f|_R$ is integrable.
	\end{proof}	
	\section{Integrals of Continuous functions}
	\par	Later we will prove a much more general result, but it is useful to start with continuous functions only and prove that continuous functions are integrable. Before we get to continuous functions, let us state the following proposition, which has a very easy proof, but it is useful to emphasize as a technique.
	\begin{lemma}
		Let $R\subset \mathbb{R}^n$ be a closed rectangle and $f:R\longrightarrow \mathbb{R}$ a bounded function. If for every $\varepsilon >0$, there exists a partition $P$ of $R$ such that 
		$$U(P,f)-L(P,f)<\varepsilon,$$
		then $f\in \mathcal{R}(R)$.
	\end{lemma}	
	\begin{proof}
		Given an $\varepsilon >0$ find $P$ as in the hypothesis. Then 
		$$\overline{\int_R}f-\underline{\int_R}f \leq U(P,f)-L(P,f)<\varepsilon.$$
		As $\displaystyle{\overline{\int_R}f\geq \underline{\int_R}f}$ and the above holds for every $\varepsilon >0$, we conclude $\displaystyle{\overline{\int_R}f= \underline{\int_R}f}$ and $f\in \mathcal{R}(R).$
	\end{proof}
	We say a rectangle $R=[\xi_1,\eta_1]\times [\xi_2,\eta_2]\times \cdots \times [\xi_n,\eta_n]$ has longest side at most $\alpha$ if $\eta_k-\xi_k\geq \alpha$ for all $k=1,2,\dots,n.$	
	\begin{lemma}
		If a rectangle $R\subset \mathbb{R}^n$ has longest side at most $\alpha$. Then for any $\xi,\eta\in R$,
		$$\|\xi-\eta\|\leq \sqrt{n}\alpha.$$
	\end{lemma}
	\begin{proof}
		\begin{align*}
			\|\xi-\eta\|=\sqrt{(x_1-y_1)^2+(x_2-y_2)^2+\cdots+(x_n-y_n)^2}\leq \sqrt{(\eta_1-\xi_1)^2+(\eta_2-\xi_2)^2+\cdots + (\eta_n-\xi_n)^2}\\
			\leq	\sqrt{\alpha^2+\alpha^2+\cdots + \alpha^2}=\sqrt{n}\alpha.
		\end{align*}
	\end{proof}
	\begin{theorem}
		Let $R\subset \mathbb{R}^n$ be a closed rectangle and $f:R\longrightarrow \mathbb{R}$ a continuous function, then $f\in \mathcal{R}(R)$.
	\end{theorem}
	\begin{proof}
		The proof is analogous to the one variable proof with some complications. The set $R$ is closed and bounded and hence compact. So $f$ is not just continuous, but in fact uniformly continuous. Let $\varepsilon >0$ be given. Find a $\delta >0$ such that $\|\xi-\eta\|<\delta$ implies $|f(\xi)-f(\eta)|<\frac{\varepsilon}{V(R)}$.\\
		\par Let $P$ be a partition of $R$ such that the longest side of any sub-rectangle is strictly less than $\frac{\delta}{\sqrt{n}}$. Then for all $\xi,\eta \in R_k$ for a sub-rectangle $R_k$ of $P$ we have, by the proposition above, $\|\xi-\eta\|<\sqrt{n}\frac{\delta}{\sqrt{n}}=\delta$.\\
		Therefore,
		$$f(\xi)-f(\eta)\leq |f(\xi)-f(\eta)|<\frac{\varepsilon}{V(R)}.$$
		As $f$ is continuous on $R_k$, it attains a maximum and a minimum on this interval. Let $\xi$ be a point where $f$ attains the maximum and $\eta$ be a point where $f$ attains the minimum. Then $f(\xi)=M_k$ and $f(\eta)=m_k$ in the notation from the definition of the integral. Therefore, 
		$$M_i-m_i=f(\xi)-f(\eta)<\frac{\varepsilon}{V(R)}.$$
		And so 
		$$U(P,f)-L(P,f)=\sum_{k=1}^N M_kV(R_k) - \sum_{k=1}^Nm_kV(R_k) = \sum_{k=1}^N(M_k-m_k)V(R_k)<\frac{\varepsilon}{V(R)}\sum_{k=1}^NV(R_k)=\varepsilon.$$
		As $\varepsilon >0$ was arbitrary,
		$$\overline{\int_a^b}f=\underline{\int_a^b}f,$$
		and $f$ is Riemann integrable on $R$.
	\end{proof}	
	\section{Integration of functions with Compact Support}
	Let $U\subset \mathbb{R}^n$ be an open set and $f:U\longrightarrow \mathbb{R}$ be a function. We say the \textit{support} of $f$ is the set 
	$$ supp(f):= \overline{\{x\in U: f(x)\neq 0\}}.$$	
	That is, the support is the closure of the set of points where the function is non-zero. The closure is in $U$, that is in particular $supp(f)\subset U$. So for a point $x\in U$ not in the support we have that $f$ is constantly zero in a whole neighborhood of $x$. 
	\par A function $f$ is said to have \textit{compact support} if $supp(f)$ is a compact set. We will mostly consider the case when $U=\mathbb{R}^n$.
	\begin{lemma}
		Suppose $f:\mathbb{R}^n \longrightarrow \mathbb{R}$ be a function with compact support. If $R$ is a closed rectangle such that $supp(f)\subset R^o$ where $R^o$ is the interior of $R$, and $f$ is integrable over $R$ then for any other closed rectangle $S$ with $supp(f)\subset S^o$, the function $f$ is integrable over $S$ and 
		$$\int_S f = \int_R f $$
	\end{lemma}
	\begin{proof}
		The intersection of closed rectangles is again a closed rectangle (or empty). Therefore we can take $\tilde{R}=R\cap S$ be the intersection of all rectangles containing $supp(f)$. If $R$ is the empty set, then $supp(f)$ is the empty set and $f$ is identically zero and the lemma becomes trivial. So suppose that $\tilde{R}$ is non-empty. As $\tilde{R}\subset R$, we know that $f$ is integrable over $\tilde{R}$. Furthermore $\tilde{R}\subset S$. Given $\varepsilon >0$, take $\tilde{P}$ to be a partition of $\tilde{R}$ such that 
		$$U(\tilde{P},f|_{\tilde{R}})-L(\tilde{P},f|_{\tilde{R}})<\varepsilon .$$
		Now add the endpoints of $S$ to $\tilde{P}$ to create a new partition $P$. Note that the sub-rectangles of $\tilde{P}$ are sub-rectangles of $P$ as well. Let $R_1,R_2,\dots,R_k$ be the sub-rectangles of $\tilde{P}$ and $R_{K+1},\dots,R_N$ the sub-rectangles. Note that since $supp(f)\subset \tilde{R}$, then for $k=K+1,\dots , N$ we have $supp(f)\cap R_k = \emptyset$. In other words $f$ is identically zero on $R_k$. Therefore in the notation used previously we have 
		\begin{align*}
			U(P,f|_S)-L(P,f|_S))=\sum_{k=1}^K(M_k-m_k)V(R_k)+\sum_{k=K+1}^N(M_k-m_k)V(R_k)=\sum_{k=1}^K(M_k-m_k)V(R_k)\\=U(\tilde{P},f|_{\tilde{R}})-L(\tilde{P},f|_{\tilde{R}})<\varepsilon .
		\end{align*}
		Similarly we have that $L(P,f|_S)=L(\tilde{P},f|_{\tilde{R}})$ and therefore
		$$ \int_S f = \int_{\tilde{R}}f. $$
		Since $\tilde{R}\subset R$ we also get $\displaystyle{\int_R f = \int_{\tilde{R}}f}$, or in other words $\displaystyle{\int_R f = \int_{S}f}$.
	\end{proof}
	\vspace{5mm}
	\section{Iterated Integrals and Fubini's Theorem}
	\par The Riemann integral in several variables is hard to compute from the definition. For one-dimensional integral we have the fundamental theorem of calculus and we can compute many integrals without having to appeal to the definition of the integral. We will rewrite a Riemann integral in several variables into several one-dimensional Riemann integrals by iterating. However, if $f:[0,1]^2 \longrightarrow \mathbb{R}$ is a Riemann integrable function, it is not immediately clear if the three expressions 
	$$\int_{[0,1]^2}f, \quad \int_0^1\int_0^1 f(x,y)\, dxdy, \quad and \quad \int_0^1\int_0^1 f(x,y)\, dydx $$	
	are equal, or if the last two are even well-defined.
	\subsection{Example}
	Define 
	$$f(x,y):=\begin{cases}
		1 &\text{if} \  x=1/2 \ and \ y\in \mathbb{Q},\\
		0 & \text{otherwise}.
	\end{cases}$$
	Then $f$ is Riemann integrable on $R:=[0,1]^2$ and $\displaystyle{\int_R f=0}$. Furthermore, $\displaystyle{\int_0^1\int_0^1 f(x,y)\,dxdy=0}$. However 
	$$\int_{0}^1 f(1/2,y)\, dy$$
	does not exist, so we cannot even write $\displaystyle{\int_0^1\int_0^1 f(x,y)\,dxdy=0}$.
	\begin{proof}
		Let us start with the integrability of $f$. We simply take the partition of $[0,1]^2$ where the partition in the $x$ direction is $\{0,1/2-\varepsilon,1/2+\varepsilon,1\}$ and in the $y$ direction $\{0,1\}$. The sub-rectangles of the partition are 
		$$R_1 :=[0,1/2-\varepsilon]\times [0,1], \quad R_2:= [1/2-\varepsilon,1/2+\varepsilon]\times [0,1], \quad R_3:=[1/2+\varepsilon,1]\times [0,1].$$
		We have $m_1=M_1=0,m_2=0,M_2=1,$ and $m_3=M_3=0.$ Therefore,
		$$L(P,f)= m_1(1/2-\varepsilon)\cdot 1 + m_2(2\varepsilon)\cdot 1 + m_3(1/2-\varepsilon)\cdot 1 = 0,$$
		and 
		$$U(P,f)= M_1(1/2-\varepsilon)\cdot 1 + M_2(2\varepsilon)\cdot 1 + M_3(1/2-\varepsilon)\cdot 1 = 2\varepsilon. $$
		The upper and lower sum are arbitrarily close and the lower sum is always zero, so the function is integrable and $\displaystyle{\int_R f=0}$.
		\vspace{1mm}
		\par For any $y$, the function that takes $x$ to $f(x,y)$ is zero except perhaps at a single point $x=1/2$. We know that such a function is integrable and $\displaystyle{\int_0^1f(x,y)\,dx=0}$. Therefore,  $\displaystyle{\int_0^1\int_0^1 f(x,y)\,dxdy=0}$. 	
		\vspace{3mm}
		\par However if $x=1/2$, the function that takes $y$ to $f(1/2,y)$ is the non-integrable function that is 1 on the rationals and 0 on the irrationals.
	\end{proof}
	\vspace{3mm}
	\begin{theorem}(Fubini version I).
		Let $R\times S\subset \mathbb{R}^n \times \mathbb{R}^m$ be a closed rectangle and \\ $f:R\times S \longrightarrow \mathbb{R}$ be integrable. The functions $g:R\longrightarrow \mathbb{R}$ and $h:R\longrightarrow \mathbb{R}$ defined by 
		$$g(x):=\underline{\int_S}f_x \quad and \quad h(x):= \overline{\int_S}f_x $$
		are integrable over $R$ and 
		$$\int_R g=\int_R h=\int_{R\times S}f.$$
	\end{theorem}	
	In other words
	$$\int_{R\times S} f = \int_R \left(\underline{\int_S}f(x,y)\,dy\right)\,dx = \int_R\left(\overline{\int_S}f(x,y)\,dy\right)\,dx.$$
	If it turns out that $f_x$ is integrable for all $x$, for example when $f$ is continuous, then we obtain the more familiar 
	$$\int_{R\times S}f = \int_R \int_S f(x,y)\,dydx.$$
	\begin{proof}
		Let $P$ be a partition of $R$ and $P'$ be a partition of $S$. Let $R_1,R_2,\dots,R_N$ be the sub-rectangles of $P$ and $R'_1,R'_2,\dots,R'_K$ be the sub-rectangles of $P'$. Then $P\times P'$ is the partition whose sub-rectangles are $R_j\times R'_k$ for all $1\leq j \leq N$ and all $1\leq k \leq K$.
		\par Let 
		$$ m_{j,k}:=\inf_{(x,y)\in R_j \times R'_k}f(x,y).$$
		We notice that $V(R_j\times R'_k)=V(R_j)V(R'_k)$ and hence 
		$$ L(P\times P',f)=\sum_{j=1}^N\sum_{k=1}^K m_{j,k}V(R_j\times R'_k)=\sum_{j=1}^N\left(\sum_{k=1}^K m_{j,k}V(R'_k)\right)V(R_j).$$
		If we let 
		$$m_k(x):=\inf_{y\in R'_k}f(x,y)=\inf_{y\in R'_k}f_x(y),$$
		then of course if $x\in R_j$ then $m_{j,k}\leq m_k(x)$. Therefore
		$$\sum_{k=1}^Km_{j,k}V(R'_k)\leq \sum_{k=1}^Km_k(x)V(R'_k)=L(P',f_x)\leq \underline{\int_S}f_x=g(x).$$
		As we have the inequality for all $x\in R_j$ we have 
		$$ \sum_{k=1}^K m_{j,k}V(R'_k)\leq \inf_{x\in R_j} g(x).$$
		We thus obtain 
		$$L(P\times P',f)\leq \sum_{j=1}^N\left(\inf_{x\in R_j}g(x)\right)V(R_j)=L(P,g).$$
		Similarly $U(P\times P',f)\geq U(P,h)$. Putting this together we have
		$$L(P\times P',f)\leq L(P,g)\leq U(P,g)\leq U(P,h)\leq U(P\times P',f).$$
		And since f is integrable, it must be that $g$ is integrable as 
		$$U(P,g)-L(P,g)\leq U(P\times P',f)-L(P\times P',f),$$
		and we can make the right hand side arbitrarily small. Furthermore as $L(P\times P',f)\leq L(P,g)\leq U(P\times P',f)$ we must have that  
		$$\int_R g = \int_{R\times S}f.$$
		Similarly we have 
		$$L(P\times P',f)\leq L(P,g) \leq L(P,h)\leq U(P,h)\leq U(P\times P',f),$$
		and hence 
		$$U(P,h)-L(P,h)\leq U(P\times P',f)-L(P\times P',f).$$
		So if $f$ is integrable so is $h$, and as $L(P\times P'.f)\leq L(P,h)\leq U(P\times P',f)$ we must have that 
		$$\int_R h = \int_{R\times S}f.$$
	\end{proof}
	\vspace{3mm}
	\begin{theorem}(Fubini Version II).
		Let $R\times S \subset \mathbb{R}^n \times \mathbb{R}^m$ be a closed rectangle and $f:R\times S \longrightarrow \mathbb{R}$ be integrable. The function $g:S\longrightarrow \mathbb{R}$ and $h:S \longrightarrow \mathbb{R}$ defined by
		$$g(y):= \underline{\int_R}f_y \quad and \quad h(y):=\overline{\int_R}f_y$$
		are integrable over $S$ and 
		$$ \int_S g=\int_S h=\int_{R\times S}f.$$
	\end{theorem}
	That is we also have 
	$$\int_{R\times S} f = \int_S \left(\underline{\int_R}f(x,y)\,dy\right)\,dx = \int_S\left(\overline{\int_R}f(x,y)\,dy\right)\,dx.$$	
	\par Next suppose that $f_x$ and $f_y$ are integrable for simplicity. For example, suppose that $f$ is continuous. Then by putting the two versions together we obtain the familiar
	$$\int_{R\times S} f=\int_R\int_S f(x,y)\,dydx = \int_S\int_Rf(x,y)\,dxdy.$$
	\par Often the Fubini's theorem is stated in two dimensions for a continuous function \\ $f:R\longrightarrow \mathbb{R}$ on a rectangle $R=[a,b]\times [c,d]$. Then the Fubini's theorem states that 
	$$\int_R f = \int_a^b\int_c^d f(x,y)\,dydx = \int_c^d\int_a^b f(x,y)\,dxdy.$$
	And the Fubini's theorem is commonly thought of as the theorem that allows us to swap the order of iterated integrals.
	\par Repeatedly applying Fubini's theorem gets us the following corollary; Let $R:=[a_1,b_1]\times [a_2,b_2]\times \cdots \times [a_n,b_n]\subset \mathbb{R}^n$ be a closed rectangle and let $f:R\longrightarrow \mathbb{R}$ be continuous. Then
	$$\int_R f = \int_{a_1}^{b_1}\int_{a_2}^{b_2}\cdots \int_{a_n}^{b_n}f(\zeta_1,\zeta_2,\zeta_3,\dots, \zeta_n)\,d\zeta_n d\zeta_{n-1}\cdots d\zeta_1.$$
	Clearly we can also switch the order of integration to any order we please. We can also relax the continuity requirement by making sure that all the intermediate functions are integrable, or by using upper or lover integrals.
	\subsection{Examples.}
	1. Compute, 
	$$\int_0^1\int_0^1\frac{x^2-y^2}{(x^2+y^2)^2}\,dxdy \quad and \quad \int_0^1\int_0^1\frac{x^2-y^2}{(x^2+y^2)^2}\,dydx. $$
	Firstly, we consider the "innermost" integral. 
	\begin{align*}
		\int_0^1\frac{x^2-y^2}{(x^2+y^2)^2}\,dy = 	\int_0^1\frac{x^2+y^2-2y^2}{(x^2+y^2)^2}\,dy= \int_0^1\frac{dy}{x^2+y^2}-\int_0^1\frac{2y^2}{(x^2+y^2)^2}\,dy\\=\int_0^1\frac{dy}{x^2+y^2}+\int_0^1y\frac{d}{dy}\left(\frac{1}{x^2+y^2}\right)=\int_0^1\frac{dy}{x^2+y^2}+\left(\left[\frac{y}{x^2+y^2}\right]^{y=1}_{y=0}-\int_0^1\frac{dy}{x^2+y^2}\right)\\ = \frac{1}{1+x^2}
	\end{align*}
	This takes care of the "innermost" integral with respect to $y$; now we do the "outermost" integral with respect to $x$ and we have 
	$$\int_0^1\frac{1}{1+x^2}\,dx = \frac{\pi}{4}.$$
	Thus, we have 
	$$\int_0^1\int_0^1\frac{x^2-y^2}{(x^2+y^2)^2}\,dxdy=-\frac{\pi}{4} \quad and \quad \int_0^1\int_0^1\frac{x^2-y^2}{(x^2+y^2)^2}\,dydx= \frac{\pi}{4} .$$
	\vspace{2mm}
	2. 	Compute, 
	$$\int_R |xy|\,dA$$
	where $R$ is the rectangle $0\leq x \leq 2, -1\leq y \leq 1$.\\
	Observe that the function $f(x,y)=|xy|$ is not really discontinuous; however, its formula in terms of the variables $x$ and $y$ depends on the sign of $xy$. Since $x$ is always positive within the rectangle $R$, we have 
	$$f(x,y)=\begin{cases}
		-xy & -1\leq y <0,\\
		xy  &  0\leq y\leq 1.
	\end{cases}$$
	Thus, by Fubini's theorem 
	\begin{align*}
		\int_R f(x,y)\,dA=\int_0^2\left(\int_{-1}^1f(x,y)\,dy\right)\,dx=\int_0^2\left(\int_{-1}^0f(x,y)\,dy+\int_0^1f(x,y)\,dy\right)\,dx\\
		= \int_0^2\left(\int_{-1}^0-xy\,dy+\int_0^1fxy\,dy\right)\,dx=\int_0^2x\,dx =2
	\end{align*}
	\section{The Set of Riemann Integrable Functions}
	\par Before we characterize all Riemann integrable functions, we need to make a slight detour. We introduce a way of measuring the size of sets in $\mathbb{R}^n$.	
	\begin{definition}
		Let $S\subset \mathbb{R}^n$ be a subset. We define the outer measure of $S$ as 
		$$m^*(S):= \inf \sum_{j=1}^{\infty}V(R_j),$$
		where the infimum is taken over all sequences $\{R_j\}$ of open rectangles such that $S\subset \bigcup_{j=1}^{\infty}R_j$. In particular, $S$ is of measure zero or a null set if $m^*(S)=0.$
	\end{definition}	
	\par We will only need measure zero sets and so we focus on these. Note that $S$ is of measure zero if and only if for every $\varepsilon>0$ there exists a sequence of open rectangles $\{R_j\}$ such that 
	$$S\subset \bigcup_{j=1}^{\infty}R_j \quad and \quad \sum_{j=1}^{\infty}V(R_j)<\varepsilon.$$
	Further-more, if $S$ is measure zero and $S'\subset S$, then $S'$ is of measure zero. We can in fact use the same exact rectangles. 
	\par We can also use balls and it is sometimes more convenient. In fact we can choose balls no bigger than a fixed radius.
	\section{Oscillation and Continuity}
	Let $S\subset \mathbb{R}^n$ be a set and $f:S\longrightarrow \mathbb{R}$ a function. Instead of just saying that $f$ is or is not continuous at a point $x\in S$ we need to able to quantify how discontinuous $f$ is as a function at $x$. For any $\delta >0$ we define the oscillation of $f$ on the $\delta$-ball in a subset topology that is $B_S(x,\delta)=B_{\mathbb{R}^n}(x,\delta)\cap S$ as
	$$o(f,x,\delta):= \sup_{y\in B_S(x,\delta)} f(y)-\inf_{y\in B_S(x,\delta)}f(y)=\sup_{y_1,y_2 \in B_S(x,\delta)}(f(y_1)-f(y_2)).$$
	That is, $o(f,x,\delta)$ is the length of the smallest interval that contains the image $f(B_S(x,\delta))$. Clearly $o(f,x,\delta)\geq 0$ and notice $o(f,x,\delta)\leq o(f,x,\delta')$ whenever $\delta<\delta'$. Therefore, the limit as $\delta \to 0$ from the right exists and we define the oscillation of a function $f$ at $x$ as 
	\vspace{-3mm}
	$$o(f,x):= \lim_{\delta\to 0^+}o(f,x,\delta)=\inf_{\delta>0}o(f,x,\delta).$$	
	\begin{lemma}
		$f:S\longrightarrow \mathbb{R}$ is continuous at $x\in S$ if and only if $o(f,x)=0$.
	\end{lemma}
	\begin{proof}
		First suppose that $f$ is continuous at $x\in S$. Then given any $\varepsilon >0$, there exists a $\delta>0$ such that for $y\in B_S(x,\delta)$ we have $|f(x)-f(y)|<\varepsilon$. Therefore if $y_1,y_2 \in B_S(x,\delta)$ then 
		$$f(y_1)-f(y_2)=f(y_1)-f(x)-(f(y_2)-f(x))<\varepsilon + \varepsilon =2\varepsilon.$$
		We take the supremum over $y_1$ and $y_2$ then, 
		$$o(f,x,\delta)=\sup_{y_1,y_2 \in B_S(x,\delta)}(f(y_1)-f(y_2))\leq 2\varepsilon.$$
		Hence, $o(f,x)=0.$
		\par On the other hand suppose that $o(x,f)=0.$ Given any $\varepsilon >0$, find a $\delta >0$ such that $o(f,x,\delta)<\varepsilon.$ If $y\in B_S(x,\delta)$ then,
		$$ |f(x)-f(y)|\leq \sup_{y_1,y_2 \in B_S(x,\delta)} (f(y_1)-f(y_2))=o(f,x,\delta)<\varepsilon.$$
	\end{proof}
	\section{The Set of Riemann Integrable Functions}	
	We have seen that continuous functions are Riemann integrable, but we also know that certain kinds of discontinuities are allowed. It turns out that as long as the discontinuities happen on a set of measure zero, the function is integrable and vice versa.
	\begin{theorem}(Riemann-Lebesgue).
		Let $R\subset \mathbb{R}^n$ be a closed rectangle and $f:R\longrightarrow \mathbb{R}$ a bounded function. Then $f$ is Riemann integrable if and only if the set of discontinuities of $f$ is of measure zero (a null set).
	\end{theorem}
	\begin{proof}
		Let $S\subset R$ be the set of discontinuities of $f$. That is $S=\{x\in R : o(f,x)>0\}$. The trick of this proof is to isolate the bad set into a small set of sub-rectangles of a partition. There are only finitely many sub-rectangles of a partition, so we will wish to use compactness. If $S$ is closed, then it would be compact and we could cover it by small rectangles as it is of measure zero. Unfortunately, in general $S$ is not closed so we need to work a little harder.
		\par For every $\varepsilon >0$, define 
		$$ S_{\varepsilon}:=\{x\in R: o(f,x)\geq \varepsilon\}.$$
		Thus, $S_{\varepsilon}$ is closed and as it is a subset of $R$ which is bounded, $S_{\varepsilon}$ is compact. Further-more, $S_{\varepsilon}\subset S$ and $S$ is of measure zero. We see that there are finitely many open rectangles $S_1,S_2,\dots, S_k$ that cover $S_{\varepsilon}$ and $\displaystyle{\sum V(S_j)<\varepsilon}.$
		\par The set $T=R\setminus (S_1\cup S_2\cup \cdots \cup S_k)$ is closed, bounded, and therefore compact. Further-more for $x\in T$, we have $o(f,x)<\varepsilon$. Hence for each $x\in T$, there exists a small closed rectangle $T_x$ with $x$ in the interior of $T_x$, such that 
		$$\sup_{y\in T_x}f(y)-\inf_{y\in T_x}f(y)<2\varepsilon.$$
		The interiors of the rectangles $T_x$ cover $T$. As $T$ is compact there exist finitely many such rectangles $T_1,T_2,\dots, T_m$ that covers $T$. 
		\par Now take all the rectangles $T_1,T_2,\dots,T_m$ and $S_1,S_2,\dots,S_k$ and construct a partition out of their endpoints. That is construct a partition $P$ with sub-rectangles $R_1,R_2,\dots,R_p$ such that every $R_j$ is contained in $T_{\ell}$ for some $\ell$ or the closure of $S_{\ell}$ for some $\ell$. Suppose we order the rectangles so that $R_1,R_2,\dots,R_q$ are those that are contained in some $T_{\ell}$, and $R_{q+1},R_{q+2},\dots,R_p$ are the rest. In particular, we have 
		$$\sum_{j=1}^q V(R_j)\leq V(R) \quad and \quad \sum_{j=q+1}^p V(R_j)\leq \varepsilon.$$
		Let $m_j$ and $M_j$ be the $\inf$ and $\sup$ over $R_j$ as before. If $R_j\subset T_{\ell}$ for some $\ell$, then $(M_j-m_j)<2\varepsilon.$ Let $B\in \mathbb{R}$ such that $|f(x)|\leq B$ for all $x\in R$, so $(M_j-m_j)<2B$ over all rectangles. Then 
		\begin{align*}
			U(P,f)-L(P,f)=\sum_{j=1}^p (M_j-m_j)V(R_j)=\left(\sum_{j=1}^q (M_j-m_j)V(R_j)\right)+\left(\sum_{j=q+1}^p(M_j-m_j)V(R_j)\right)\\
			\leq \left(\sum_{j=1}^q 2\varepsilon V(R_j)\right)+\left(\sum_{j=q+1}^p 2BV(R_j)\right)\leq 2\varepsilon V(R)+2B\varepsilon = \varepsilon(2V(R)+2B).
		\end{align*}
		Clearly, we can make the right hand side as small as we want and hence $f$ is integrable. 
		\par For the other direction, suppose that $f$ is Riemann integrable over $R$. Let $S$ be the set of discontinuities again and now let 
		$$S_k:=\{x\in R : o(f,x)\geq 1/k\}.$$	
		Fix a $k\in \mathbb{N}$. Given an $\varepsilon >0$, find a partition $P$ with sub-rectangles $R_1,R_2,R_3,\dots,R_p$ such that 
		$$U(P,f)-L(P,f)=\sum_{j=1}^p(M_j-m_j)V(R_j)<\varepsilon$$
		Suppose that $R_1,R_2,\dots,R_p$ are ordered so that the interiors of $R_1,R_2,\dots, R_q$ intersect $S_k$, while the interiors of $R_{q+1},R_{q+2},\dots, R_p$ are disjoint from $S_k$. If $x\in R_j\cap S_k$ and $x$ is in the interior of $R_j$ so sufficiently small balls are completely inside $R_j$, then by definition of $S_k$ we have $M_j-m_j \leq 1/k$. Then 
		$$\varepsilon > \sum_{j=1}^p (M_j-m_j)V(R_j)\geq \sum_{j=1}^q (M_j-m_j)V(R_j)\leq \frac{1}{k}\sum_{j=1}^q V(R_j)$$
		In other words $\displaystyle{\sum_{j=1}^q V(R_j)<k\varepsilon}$. Let $G$ be the set of all boundaries of all the sub-rectangles of $P$. The set $G$ is of measure zero. Let $R^o_j$ denote the interior of $R_j$, then 
		$$S_k \subset R^o_1 \cup R^o_2 \cup R^o_3 \cup \cdots \cup R^o_q \cup G.$$
		As $G$ can be covered by open rectangles arbitrarily small volume, $S_k$ must be of measure zero. As 
		$$S= \bigcup_{k=1}^{\infty}S_k$$
		and a countable union of measure zero sets is of measure zero. $S$ is of measure zero.
	\end{proof}	
	\section{Application to Summation of Series \\ in $\mathbb{R}$}	
	\par The theory of limits of finite approximations was made precise by the German mathematician \textit{Bernhard Riemann}. We now introduce the notion of a Riemann sum, which underlies the theory of definite integrals.
	We begin with an arbitrary bounded function $f$ defined on a closed interval $[a,b]$. $f$ may have negative as well as positive values. We subdivide the interval [a,b] into subintervals, not necessarily of equals widths or lengths, and form sums in the same way as for the finite approximations.
	\par In each subinterval we select some point. The point chosen in the $kth$ subinterval $[x_{k-1},x_k]$ is called $c_k$. Then on each subinterval we stand a vertical rectangle that stretches from the $x-axis$ to touch the curve at $(c_k,f(c_k))$. These rectangles can be above or below the $x-axis$, depending on whether $f(c_k)$ is positive or negative, or on the $x-axis$ if $f(c_k)=0$ .
	On each subinterval we form the product $f(c_k)\cdot \Delta x_k$. This product is positive, negative, or zero, depending on the sign of $f(c_k)$.
	\par Finally we sum all these products to get 
	$$S_p = \sum_{k=1}^{n}f(c_k)\Delta x_k.$$	
	The sum $S_n$ is called a \textbf{Riemann sum for f on the interval [a,b]}. There are many such sums, depending on the partition $P$ we choose, and the choices of the points $c_k$ in the subintervals.
	Here we choose $\Delta x = (b-a)/n$ to partition $[a,b]$, and then choose the point $c_k$ to be the right-handed endpoint of each subinterval when forming the Riemann sum. This choice leads to the Riemann sum formula 
	$$S_n = \sum_{k=1}^n f\left(a+k\frac{(b-a)}{n}\right)\cdot \left(\frac{b-a}{n}\right)$$
	\begin{definition}
		Let $f(x)$ be a function defined on a closed interval $[a,b]$. We say that a number $J$ is the definite integral of $f$ over $[a,b]$ and that $J$ is the limit of the Riemann sum $\displaystyle{\sum_{k=1}^n f(c_k)\Delta x_k}$ if the following condition is satisfied: 
		\par Given any number $\varepsilon >0$ there is a corresponding number $\delta >0$ such that for every partition $P$ of $[a,b]$ with $\|P\|<\delta$ and any choice of $c_k$ in $[x_{k-1},x_k]$, we have 
		$$\left| \ \sum_{k=1}^n f(c_k) \Delta x_k - J \  \right |<\varepsilon .$$
	\end{definition}
	The symbol for the number $J$ in the definition of the definite integral is 
	$$\int_a^b f(x)\,dx.$$
	So we can write, 
	$$J= \int_a^b f(x)\, dx = \lim_{n\to \infty}\sum_{k=1}^n f(c_k)\left(\frac{b-a}{n}\right)=\lim_{n\to \infty}\sum_{k=1}^n f\left(a+k\frac{(b-a)}{n}\right)\left(\frac{b-a}{n}\right)	$$
	\section{Examples.}
	1. Evaluate 
	$$\lim_{n\to \infty} \left(\frac{1}{n}\left(\frac{1}{n}\right)^2+\frac{1}{n}\left(\frac{2}{n}\right)^2+\frac{1}{n}\left(\frac{3}{n}\right)^2+\cdots +\frac{1}{n}\left(\frac{n}{n}\right)^2\right)$$
	Now, 
	\begin{align*}
		\lim_{n\to \infty} \left(\frac{1}{n}\left(\frac{1}{n}\right)^2+\frac{1}{n}\left(\frac{2}{n}\right)^2+\frac{1}{n}\left(\frac{3}{n}\right)^2+\cdots +\frac{1}{n}\left(\frac{n}{n}\right)^2\right)=\lim_{n\to \infty}\sum_{i=1}^n \frac{1}{n}\left(\frac{i}{n}\right)^2 \\
		=	\lim_{n\to \infty}\sum_{i=1}^n \frac{1-0}{n}\left(0+i\frac{(1-0)}{n}\right)^2 = \int_0^1 x^2\,dx = \frac{1}{3}
	\end{align*}
	That is 
	$$\lim_{n\to \infty} \left(\frac{1}{n}\left(\frac{1}{n}\right)^2+\frac{1}{n}\left(\frac{2}{n}\right)^2+\frac{1}{n}\left(\frac{3}{n}\right)^2+\cdots +\frac{1}{n}\left(\frac{n}{n}\right)^2\right)=\frac{1}{3}$$	
	2. Evaluate 
	$$\lim_{n\to \infty}\left(\frac{3}{n}\left(2+\frac{3}{n}\right)^2+\frac{3}{n}\left(2+\frac{6}{n}\right)^2+\frac{3}{n}\left(2+\frac{9}{n}\right)^2+\cdots + \frac{3}{n}\left(2+\frac{3n}{n}\right)^2\right)	$$
	Notice,
	\begin{align*}
		\lim_{n\to \infty}\left(\frac{3}{n}\left(2+\frac{3}{n}\right)^2+\frac{3}{n}\left(2+\frac{6}{n}\right)^2+\frac{3}{n}\left(2+\frac{9}{n}\right)^2+\cdots + \frac{3}{n}\left(2+\frac{3n}{n}\right)^2\right)= \lim_{n\to \infty} \sum_{i=1}^n \frac{3}{n}\left(2+\frac{3i}{n}\right)^2\\
		= \lim_{n\to \infty} \sum_{i=1}^n \frac{(5-2)}{n}\left(2+i\frac{(5-2)}{n}\right)^2=\int_2^5 x^2\,dx = 39
	\end{align*}	
	It follows immediately that,
	$$\lim_{n\to \infty}\left(\frac{3}{n}\left(2+\frac{3}{n}\right)^2+\frac{3}{n}\left(2+\frac{6}{n}\right)^2+\frac{3}{n}\left(2+\frac{9}{n}\right)^2+\cdots + \frac{3}{n}\left(2+\frac{3n}{n}\right)^2\right)=39	$$
	3. Find a closed form for,
	$$\lim_{n\to \infty}\left(\frac{b-a}{n\left(a+\frac{b-a}{n}\right)^2}+\frac{b-a}{n\left(a+2\frac{b-a}{n}\right)^2}+\frac{b-a}{n\left(a+3\frac{b-a}{n}\right)^2}+\cdots + \frac{b-a}{n\left(a+n\frac{b-a}{n}\right)^2}\right)$$
	Observe that, 
	$$\lim_{n\to \infty}\left(\frac{b-a}{n\left(a+\frac{b-a}{n}\right)^2}+\frac{b-a}{n\left(a+2\frac{b-a}{n}\right)^2}+\frac{b-a}{n\left(a+3\frac{b-a}{n}\right)^2}+\cdots + \frac{b-a}{n\left(a+n\frac{b-a}{n}\right)^2}\right)$$ $$ =\lim_{n\to \infty}\frac{b-a}{n}\sum_{k=1}^n\left(a+\frac{b-a}{n}k\right)^{-2}$$
	Now let 
	$$S_n = 	\frac{b-a}{n}\sum_{k=1}^n\left(a+\frac{b-a}{n}k\right)^{-2}$$
	Observe that 
	$$\frac{b-a}{n}\sum_{k=1}^n\left(a+\frac{b-a}{n}k\right)^{-1}\left(a+\frac{b-a}{n}(k+1)\right)^{-1}\leq S_n $$ 
	Thus, $$S_n\leq \frac{b-a}{n}\sum_{k=1}^n\left(a+\frac{b-a}{n}k\right)^{-1}\left(a+\frac{b-a}{n}(k-1)\right)^{-1} $$
	And decomposing into partial fractions we have, 
	\begin{align*}
		\sum_{k=1}^n \left\{\left(a+\frac{b-a}{n}k\right)^{-1}-\left(a+\frac{b-a}{n}(k+1)\right)^{-1}\right\}\leq S_n 
	\end{align*}	
	Likewise 
	\begin{align*}
		S_n \leq 	\sum_{k=1}^n \left\{\left(a+\frac{b-a}{n}(k-1)\right)^{-1}-\left(a+\frac{b-a}{n}k\right)^{-1}\right\}
	\end{align*}	
	Notice the sums are telescoping, then it follows that 
	\begin{align*}
		\left(a+\frac{b-a}{n}\right)^{-1}-\left(a+\frac{b-a}{n}(n+1)\right)^{-1}\leq S_n \leq \frac{1}{a}-\frac{1}{b}
	\end{align*}	
	By the squeeze theorem, we get that 
	$$\lim_{n\to \infty} S_n = \lim_{n\to \infty}\left(\frac{b-a}{n\left(a+\frac{b-a}{n}\right)^2}+\frac{b-a}{n\left(a+2\frac{b-a}{n}\right)^2}+\frac{b-a}{n\left(a+3\frac{b-a}{n}\right)^2}+\cdots + \frac{b-a}{n\left(a+n\frac{b-a}{n}\right)^2}\right)=\frac{1}{a}-\frac{1}{b}$$
	But we can make things more easier simply by noting the Riemann sum
	$$\int_a^b x^{-2}\,dx = \lim_{n\to \infty}\left(\frac{b-a}{n\left(a+\frac{b-a}{n}\right)^2}+\frac{b-a}{n\left(a+2\frac{b-a}{n}\right)^2}+\frac{b-a}{n\left(a+3\frac{b-a}{n}\right)^2}+\cdots + \frac{b-a}{n\left(a+n\frac{b-a}{n}\right)^2}\right)$$
	
	4. Show that 
	$$\int_0^{\pi}\log(a^2+b^2-2ab\cos(x))\,dx=2\pi \log(\max\{a,b\}) \quad (\text{Bronstein Integral})$$
	Now, suppose WLOG $b>a$. Then with $c=b/a$ we have 
	\begin{align*}
		I=\int_0^{\pi}\log(a^2+b^2-2ab\cos(x))\,dx = \int_{0}^{\pi}\log(a^2)\,dx + \int_0^{\pi}\log(1+(b/a)^2-2(b/a)\cos(x))\,dx\\
		= 2\pi\log(a)+\underbrace{\int_0^{\pi}\log(1+c^2-2c\cos(x))\,dx}_{J}
	\end{align*}	
	Now we can evaluate the second integral on the RHS as the limit of a \textit{Riemann sum}	
	\begin{align*}
		J=\int_0^{\pi}\log(1+c^2-2c\cos(x))\,dx = \lim_{n\to \infty}\frac{\pi}{n}\sum_{j=1}^n\log\left(1+c^2-2c\cos\left(\frac{\pi(j-1)}{n}\right)\right)\\
		= \lim_{n\to \infty}\frac{\pi}{n}\log\left(\prod_{j=1}^n\left(1+c^2-2c\cos\left(\frac{\pi(j-1)}{n}\right)\right)\right)=\lim_{n\to \infty}\log((1-c)^2)\\+\lim_{n\to \infty}\frac{\pi}{n}\log\left(\prod_{j=2}^n\left(1+c^2-2c\cos\left(\frac{\pi(j-1)}{n}\right)\right)\right)=\lim_{n\to \infty}\frac{\pi}{n}\log\left(\prod_{j=1}^{n-1}\left(1+c^2-2c\cos\left(\frac{j\pi}{n}\right)\right)\right)
	\end{align*}	
	Notice we can factor, 
	$$1+c^2-2c\cos\left(\frac{j\pi}{n}\right)= [c-\exp(i\pi j/n)][c-\exp(-i\pi j/n)],$$	
	Hence, recalling that $c=b/a>1$, it follows that 
	\begin{align*}
		J=\lim_{n\to \infty}\frac{\pi}{n}\log\left(\frac{c^{2n}-1}{c^2-1}\right)=\pi \lim_{n\to \infty}\log\left(\frac{c^{2n}-1}{c^2-1}\right)^{1/n}=\pi \lim_{n\to \infty}\log\left[c^2\left(\frac{c^{2n}-1}{c^2-1}\right)^{1/n}\right]\\
		= \pi \log(c^2) + \log\left[\lim_{n\to \infty}\left(\frac{c^{2n}-1}{c^2-1}\right)^{1/n}\right]=\pi \log(c^2)
	\end{align*}
	Thus, 
	$$I=2\pi \log(a)+\pi\log(c^2)=2\pi\log(a)+\pi\log((b/a)^2)=2\pi \log(b)$$
	Hence, we conclude 
	$$\int_0^{\pi}\log(a^2+b^2-2ab\cos(x))\,dx=2\pi \log(\max\{a,b\})$$

	\footnote{\emph{Riemann : If only I had the theorems! Then I should find the proofs easily enough.}}

\end{document}